\documentclass[review]{elsarticle}

\usepackage{lineno,hyperref}
\modulolinenumbers[5]

\journal{  }

\usepackage{graphicx}
\usepackage{subfigure}
\usepackage{dsfont} 
\usepackage{amssymb} 
\usepackage{mathrsfs} 
\usepackage{amsmath}

\begin{document}

\begin{frontmatter}

\title{Maximal Likely Phase Lines for a Reduced Ice Growth Model\tnoteref{mytitlenote}}
\tnotetext[mytitlenote]{This work was partly supported by the NSF grant 1620449,  and NSFC grants 11531006 and 11771449.}

\author[mymainaddress]{Athanasios Tsiairis}
\ead{thtsiairis@hotmail.com}

\author[mymainaddress]{Pingyuan Wei}
\ead{weipingyuan@hust.edu.cn}

\author[mymainaddress]{Ying Chao\corref{mycorrespondingauthor}}
\ead{yingchao1993@hust.edu.cn}

\author[mymainaddress,mysecondaryaddress]{Jinqiao Duan}
\cortext[mycorrespondingauthor]{Corresponding author}
\ead{duan@iit.edu}

\address[mymainaddress]{School of Mathematics and Statistics, \& Center for Mathematical Sciences,  Huazhong University of Sciences and Technology, Wuhan 430074,  China}
\address[mysecondaryaddress]{Department of Applied Mathematics, Illinois Institute of Technology, Chicago, IL 60616, USA}

\begin{abstract}
We study the impact of Brownian noise on transitions between metastable equilibrium states in a stochastic ice sheet model. Two methods to accomplish different objectives are employed. The maximal likely trajectory by maximizing the probability density function and numerically solving the Fokker-Planck equation shows how the system will evolve over time. We have especially studied the maximal likely trajectories starting near the ice-free metastable state, and examined whether they evolve to or near the  ice-covered metastable state for certain parameters, in order to gain insights into how the ice sheet formed. Furthermore, for the transition from ice-covered metastable state to the ice-free metastable state, we study the most probable path for various noise parameters via the Onsager-Machlup least action principle. This enables us to predict and visualize the melting process of the ice sheet if such a rare event ever does take place.
\end{abstract}


\end{frontmatter}

\section{Introduction}
\label{intro}
Tremendous variety of nonlinear complex dynamical systems are subject to noisy perturbations \cite{Ar,Duan}. These noises usually play a pivotal role on setting up the dynamical behavior of the system. It is incredibly important to study the influence of noise on ice sheets as it is closely connected with the livelihood of human beings \cite{Agarwal2018,Thorndike2000}. Indeed, with the increasing global mean temperatures over the last decades, the behaviour of ice sheets, such as those on Greenland and Antarctica, become a major topic in climate research \cite{Mulder2018,Toppaladoddi2017}. Research studies show that the melting of the Greenland Glaciers has caused the increase of the sea level by 0.5 mm per year in the period 2003-2008 \cite{Woot2008} while if the Antarctic glaciers were ever going to melt the sea level would rise by an astonishing 58 meters \cite{Fretwell2013}. All these make the study of the ice sheet more imminent. More results of interactions among atmosphere, ice, land and ocean for feedbacks that are relevant to understandings of climate change could be found in  \cite{Ar,Franzke2017}. 

In this present paper, we consider an ice sheet model for the development of ice sheets with boundary on the polar sea (such as the Arctic Ocean for the Greenland ice sheet), developed by Weertman \cite{Weertman1964,Weertman1976}, with the accumulation and ablation perturbed by Brownian fluctuations. In the direction of deterministic analysis, Weertman's idea has been further extended and improved through both conceptual and numerical models; see K\"all\'en et al. \cite{Kallen1979} and Oerlemans-Van der Veen \cite{Oerlemans1984}. As noise is present, the stochastic model has two metastable states in certain parameter ranges (see Section 2), which are referred as ``ice-covered" state and ``ice-free" state respectively. Random fluctuations may lead to switching between these two states, and such transitions occur widely in not only climate model, but also biological, chemical, physical, and other systems \cite{Chen2018,Duan,OM}. In order to gain insights into the evolution trajectory of the formation and even melting of ice sheet, it is of interest to study the transitions between these two metastable states. The objective of this present paper is to study the maximal likely trajectories and most probable transition paths for such a stochastic system as time goes on. This offers the following information for the ice sheet: (i) The evolution trajectories from ice-free state, indicating the forming routes of this ice sheet. (ii) The Transition paths from ice-covered state to ice-free state, predicting the melting routes on the assumption that the ice sheet will go completely in certain time interval.

We remark that ``maximal likely trajectory" and ``most probable path" here are two mathematical terms that have essential distinction \cite{Cheng2019,OM}: Firstly, The maximal likely trajectory is an evolution trajectory starting from one initial state with unknown final state, while the most probable path is a continuous transition trajectory from one metastable state to another metastable state. Secondly, the former is determined by maximizing the probability density function at every time instant, while the latter should be understood as the probability maximizer where sample solution paths lie within a tube. Thirdly, the former is obtained via numerically solving an initial value problem, while the latter is calculated by solving the two-point boundary value problem. 

This paper is organized as follows. In Section 2, we introduce the stochastic ice sheet model influenced by Brownian motion and present the solutions of the deterministic counterpart. In Section 3, we present the methods used with some preliminary results. We further consider the most probable pathways by simulating exact model solutions under noise and compare these to the maximal likely trajectories. In Section 4, we present final results with varying parameters and noise intensity while also examine examine the transition from one state to another. More specifically we investigate the possibility of a transition from the ice-covered to the ice-free state. Finally, we summarize the above results in Section 5.

\section{Model} 
We first introduce the height-mass balance feedback in an ice sheet system influenced by the Brownian motion, and then show most probable transition pathways and maximal likely trajectories for such a system.

The development of ice sheets is governed by nonlinear processes \cite{Ar}. Consider an idealised ice sheet with length $X$.  It is natural to choose the coordinates such that the point $x = 0$ corresponds to the boundary with the polar sea; see K\"all\'en et al. \cite{Kallen1979} and Oerlemans-Van der Veen \cite{Oerlemans1984}. Indeed, we can refer to the Arctic Ocean for the Greenland ice sheet as a practical examples. 
Note that the ice can be treated as a perfectly plastic material, by horizontal stress balance \cite[Page 277-278]{Ar} in the ice, the height (or thickness) of the ice sheet $h$ satisfies
\begin{equation}\label{Eqn-ice thickness}
h(x,t)=\sqrt{\sigma}\Big(\frac{X(t)}{2}-\big|x-\frac{X(t)}{2}\big| \Big)^{\frac{1}{2}},
\end{equation}
where $\sigma$ is a yield stress parameter. The maximum height of the ice sheet is indicated by $H(t)$ and given by $H(t)=h(X/2,t)=\sqrt{\sigma X(t)/2}$.
\par
On the other hand, by the mass balance for the ice cap, we have following continuity equation
\begin{equation}\label{Eqn-continuity}
\rho_i\frac{\partial h}{\partial t}=P_i-M,
\end{equation}
where the right-hand side is the mass balance for the ice cap (the difference between accumulation $P_i$ and ablation $M$). Indeed, the mass balance depends on the distance from the polar ocean, represented by $r$, and the height of the ice sheet. Note that the ice sheet is located on the north of polar sea, we always consider $r\leqslant 0$. Some experiments also indicated that both accumulation and ablation are influenced by random fluctuations arising from the complex environment in actual situation \cite{Ar,Franzke2017}. We consider a linear relation of the form
\begin{equation}\label{Eqn-mass balance}
P_i(x,t)-M(x,t)=G(x,t)=\rho_i\beta\Big(h(x,t)-\lambda(x-r)+\varepsilon_0{\xi}(t)\Big),
\end{equation}
where $\beta>0$, $\lambda>0$ are constants, ${\xi}(t)$ is a stochastic noise, $\varepsilon_0$ is the strength of the noise. 

In this paper, we will consider stochastic noise ${\xi}(t)=\dot{B}_t$ as a Gaussian white noise, which is a special stationary stochastic process, with mean $\mathbb{E}B_t=0$ and covariance $\mathbb{E}(B_tB_s)=\delta(t-s)$, and, formally, can be understood as the ``time derivative" of Brownian motion (also called Wiener process) \cite[Page 51]{Duan}. 
\begin{table}
\caption{Parameters of the conceptual ice sheet model \cite{Ar}}
\begin{tabular}{llll}
\hline\noalign{\smallskip}
 Parameter & Meaning & Value & Unit\\
 \noalign{\smallskip}\hline\noalign{\smallskip}
   $\rho$ & ice density & 0.9 & $kg$ $m^{-3}$\\
   $\sigma$ & yield stress parameter & 6.25 & $m$\\
   $\beta$ & coefficient in the linear mass balance & $10^{-3}$ & $yr^{-1}$\\
\noalign{\smallskip}\hline
\end{tabular}
\label{tbl:table1}
\end{table}
\par
Assume that the snow accumulating on the northern half of the ice sheet flows into the Arctic Ocean or melts close to it. The evolution of the ice sheet is then governed by the mass balance conditions on the southern half of the ice sheet.
\par
Substituting $G(x,t)$ into the continuity equation (\ref{Eqn-continuity}) and integrating over the southern part of the ice sheet, we obtain (with $h=\sqrt{\sigma} (X-x)^{\frac{1}{2}}$ )
\begin{equation}
\int_{X/2}^{X}\frac{\partial h}{\partial t}dx=\sqrt{\frac{\sigma X}{2}}\frac{dX}{dt}=\beta\int_{X/2}^{X}\Big[h(x,t)-\lambda(x-r)+\varepsilon_0{\xi}(t)\Big]dx.
\end{equation}
Hence the length $X$ of the ice sheet can be determined by the following stochastic differential equation (SDE)
\begin{equation}\label{Eqn-L}
\dot{X}=f(X)+\varepsilon g(X)\xi(t),
\end{equation}
where $f(X)
=-{\frac{\beta\lambda}{\sqrt{2\sigma}}} \Big(\frac{3}{4}X^{\frac{3}{2}}-r{X}^{\frac{1}{2}} \Big)+\frac{1}{3}\beta X$, $g(X)={X}^{\frac{1}{2}}$ and $\varepsilon= \frac{\beta\varepsilon_0}{\sqrt{2\sigma}}$.

\par
\begin{figure}
  \centering     
  \subfigure[]{  
  \includegraphics[width=2.2in]{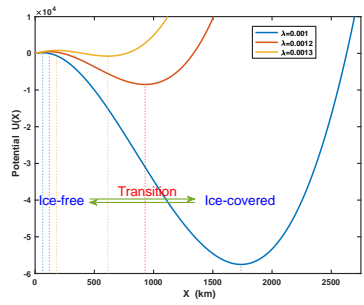}  
  }
  \subfigure[]{  
  \includegraphics[width=2.2in]{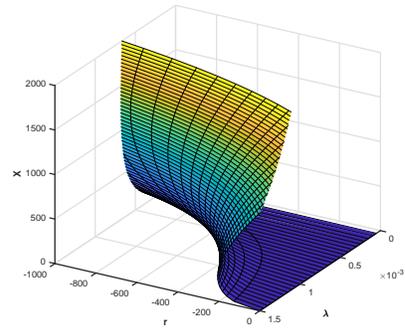}  
  }
	\caption{(Color online) (a) The potential $U$ as a function of the length $X$ with different parameter $\lambda$; (b) The cusp catastrophe surface for the deterministic ice sheet model in (\ref{Eqn-D}) with $r$ and $\lambda$ as control parameters. Indeed, equation (\ref{Eqn-D}) always has an equilibrium state, 0; when $r=0$, there exists only one additional equilibrium state, $\frac{32\sigma}{81\lambda^2} $ ($\approx2469.1 km$, if $\lambda=0.001$); but for $r<0$, there exists two additional equilibrium state,
$
X_{\pm}=\frac{4}{3}\frac{|r|}{1\mp\sqrt{\Delta}},
$
 if the discriminant $\Delta:=1+\frac{27r}{2\sigma}\lambda^2>0$. Note that $0<\Delta<1$, both of $X_{\pm}$ are bigger than 0. Thus, the state $X_{-}$ is an unstable node and the state $X_{+}$ is a stable node. 
 For $r=-250km$, the  reasonable value of $\lambda$ should be small about $0.0014$. The deterministic model with $\lambda=0.001$ has been studied in \cite{Ar}, in which these two states are $X_{-}\approx 63.9 km$ and $X_{+}\approx 1738.6 km$.  } 
	\label{fig_potential}   
\end{figure}

For the deterministic counterpart
\begin{equation}\label{Eqn-D}
\dot{X}=f(X),
\end{equation}
the vector field $f(X)$ can be rewritten as $-U^{\prime}(X)$, with the potential function $U(X)=\frac{\beta\lambda}{\sqrt{2\sigma}} \Big(\frac{3}{10}X^{\frac{5}{2}}-\frac{2}{3}r{X}^{\frac{3}{2}} \Big)-\frac{1}{6}\beta X^2$; see Fig. \ref{fig_potential}(a).
Since there are two control parameters, $\lambda$ and $r$, it is desirable to consider beyond codimension-one bifurcation events. The parameter space is now represented jointly by the $(r,\lambda)$-plane. A greater perspective to the fold bifurcation is shown in Fig. \ref{fig_potential}(b) in terms of the cusp catastrophe surface. This cusp codimension-two phenomenon is observed because we have independent control of both $r$ and $\lambda$.
\par
Without noise, in certain parameter ranges, equation (\ref{Eqn-D}) has two stable equilibrium states separated by an unstable equilibrium (which is a saddle node). In one of the stable equilibrium states, the length of ice sheet is 0 and corresponds physically to the ice sheet being melted completely. This state is usually referred to as the ice-free state. In the other equilibrium state the length of ice sheet is large and corresponds to the ice sheet being formed. This state is called ice-covered state. For different initial value of the length of ice sheet, the trajectory of $X$ may lie either in the domain of attraction of the ice-free state or in the domain of attraction of the ice-covered state. That is, the ice-free stable state 0 and ice-covered stable state $X_{+}$ are resilient (see Fig. \ref{fig_potential}(a)): the ice length states will locally be attracted to 0 or $X_{+}$, as time increases for the deterministic system. When noise is present, the system (\ref{Eqn-L}) may show switches between state 0 and state $X_{+}$, and then these two states are called metastable states. By passing through the unstable saddle state $X_{-}$, the ice length starting near the ice-free state 0 in interval $(0,X_{-})$ arrives at an ice-covered state (near $X_{+}$). 
\par
We now examine these system trajectories or orbits for the stochastic ice sheet model:
(i) How does the system evolve from ice-free situation (near $0$) to ice-covered situation (near $X_+$)? It means that we can try to understand how an ice sheet is formed.
(ii) If the ice sheet melted, how would the system transit from ice-covered situation to ice-free situation? We should also wonder how likely is such a transition for this model?

\renewcommand{\theequation}{\thesection.\arabic{equation}}
\setcounter{equation}{0}

\section{Methods}
Two different methods will be used for te purposes of this paper. I this part they will be reviewed and explained.
In the first subsection, we consider the ice sheet system (\ref{Eqn-L}) with Brownian noise (that is, $\xi(t)=\dot{B}_t$), and use the \emph{maximal likely trajectories} to discuss the evolution of trajectory starting from different initial states. In the second subsection, we focus on transition paths connecting two special fixed states for this system. The so called \emph{most probable path} will be proposed via Onsager-Machlup's function. 

\subsection{Maximal Likely Trajectory Based on Fokker-Planck Equation}
We consider the maximal likely trajectory \cite{Duan,Cheng2016} for the stochastic system (\ref{Eqn-L}), starting at an initial state $X_0$. This is reminiscent of studying a deterministic dynamical system by examining the evolution of its trajectory starting from an initial state. Each sample solution path starting at this initial state is a possible outcome of the solution path $X_t$. Then an interesting question to raise is: What is the maximal likely trajectory of $X_t$? To answer this question, we need to decide on the maximal likely position $X_{ml}$ of the system (starting at the initial point $X_0$), at every given future time $t$, and this would be the maximizer for the probability density function $p(X, t)=p(X, t; X_0, 0)$ of solution $L_t$. Indeed, the probability density function $p(X, t)$ is a surface in the $(X, t, p)$-space. At a given time instant $t$, the maximizer $X_{ml}(t)$ for $p(X, t)$ indicates the maximal likely location of this orbit at time $t$. Therefore, $X_{ml}(t)$ follows the top ridge or plateau of the surface in the $(X, t, p)$-space as time goes on, and the trajectory (or orbit) traced out by $X_{ml}(t)$ is called the \emph{maximal likely trajectory} starting at $X_0$. The maximal likely trajectories are also called `\emph{paths of mode}' in climate dynamics and data assimilation \cite{Miller1999,Cheng2019}.
\par
For the stochastic ice sheet system with Brownian noise
\begin{equation}\label{Brownian}
d{X}=f(X)dt+\varepsilon g(X)dB(t),
\end{equation}
the Fokker-Planck equation \cite{Duan} in terms of the probability density function $p(X,t)$ for the solution process $X_t$ given initial condition $X_0$ is 
\begin{equation}\label{local_pde}
 p_t = -(f(x)p(x,t))_x + \frac{1}{2}\varepsilon^2(g(x))^2p(x,t)_{xx},    \qquad p(x,0) = \delta(x-x_0).
\end{equation}
By numerically solving the Fokker-Planck equation, we can find the maximal likely position $X_{ml}(t)$ as the maximizer of $p(X, t)$ at every given time $t$. \par
\begin{figure} 
  \subfigure[]{  
  \includegraphics[width=2.3in]{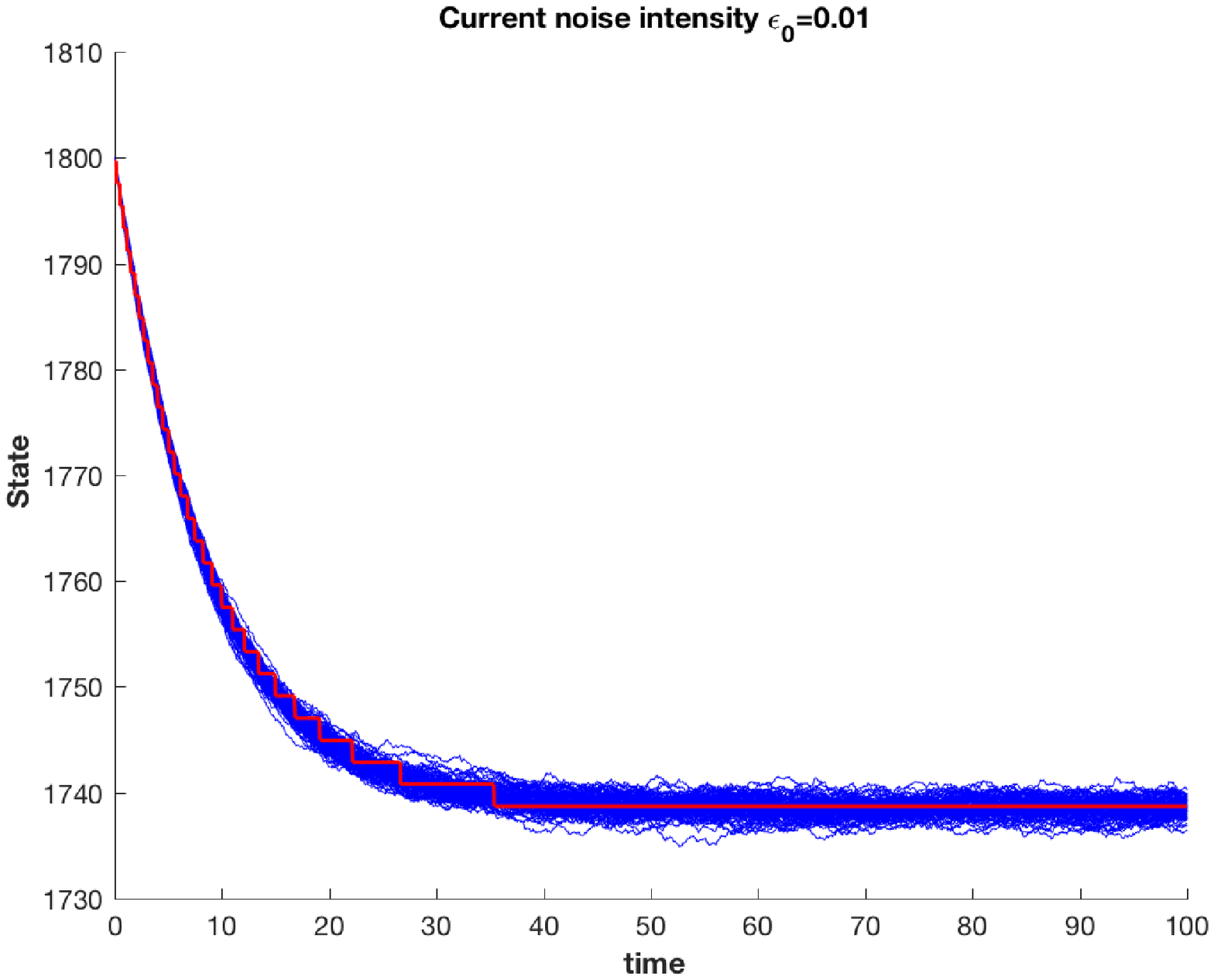}  
  }
  \subfigure[]{  
  \includegraphics[width=2.3in]{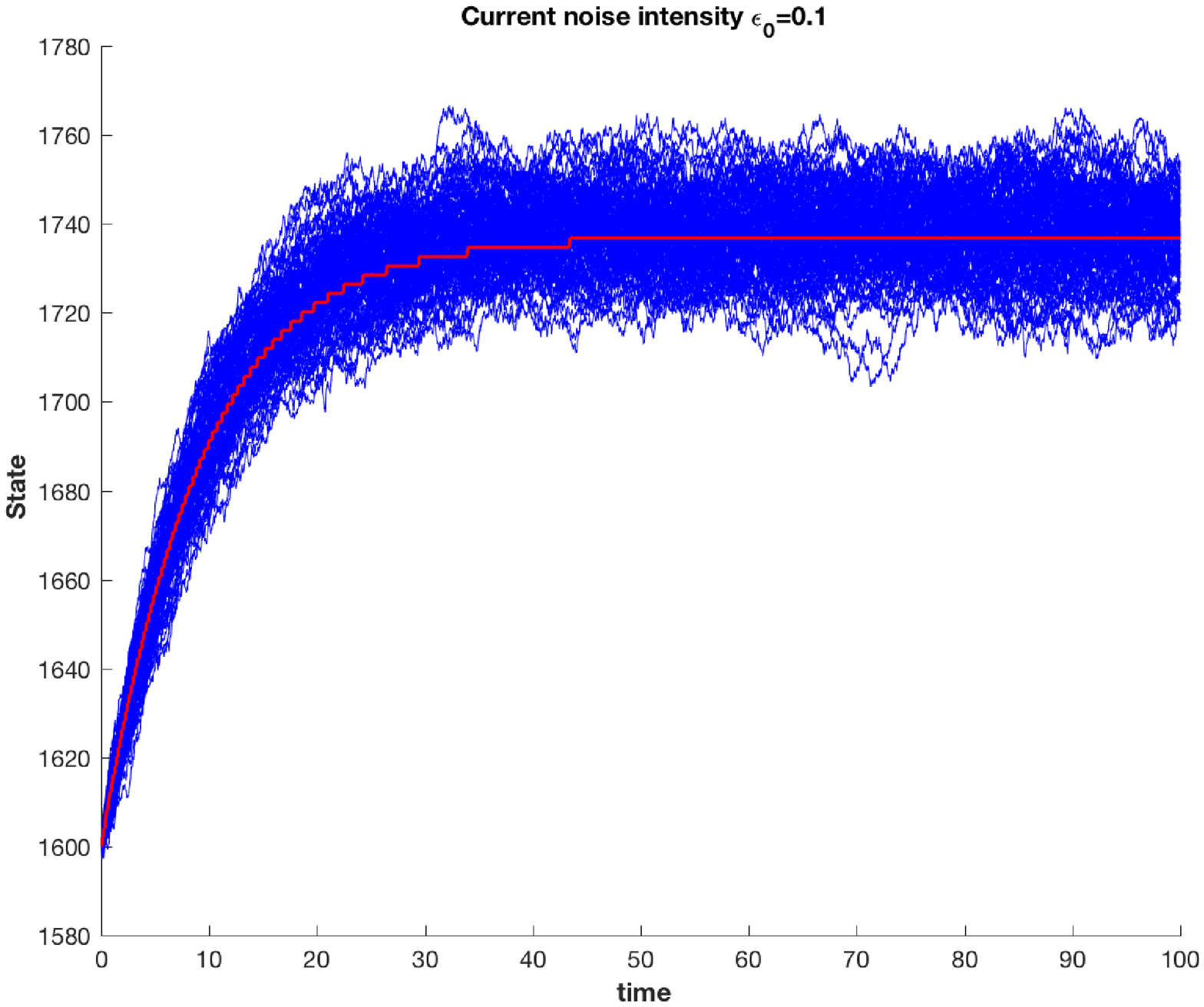}  
  }
  \caption{(Color online) Let $r=-250(km)$ and $\lambda=0.001$. Simulations of system (\ref{Brownian}) (blue lines) and the corresponding maximal likely trajectory (red line) : (a) Initial state $X_0=1800$, noise intensity $\varepsilon_0=0.01$; (b) Initial state $X_0=1600$, noise intensity $\varepsilon_0=0.1$.
  }   
  \label{fig:Actual} 
\end{figure}

We consider the unit of time as ``$kyr$" and the unit of length as ``$km$" in the numerical calculation throughout this paper. In Fig. \ref{fig:Actual}(a)(b), from around 100 stochastic simulations of system (\ref{Brownian}) with $\varepsilon_0=0.01$ and $\varepsilon_0=0.1$, respectively, we observe The maximal likely trajectory (red line) is actually located in the middle of the blue area created from all the possible outcomes. that simulation orbits are more likely to concentrate around the maximal likely trajectory $X_{ml}(t)$ (red line). This result indicates that the orbit according to the Fokker-Planck equation is actually the maximal likely evolution trajectory of the system (\ref{Brownian}).
\par
The state which attracts (or repels) all nearby orbits is referred as a \emph{maximal likely stable (unstable) equilibrium state}  \cite{Wang}, which depends on noise intensity $\varepsilon_0$ as well as the ice sheet system parameters $\lambda$, $r$. Fig. \ref{fig:Actual} shows that the maximal likely equilibrium states in these two cases are $1736.8km$ and $1734.7km$, respectively. 
We will exhibit the number and value of maximal likely stable equilibrium states for the stochastic ice sheet model more thoroughly in Section 4.

\subsection{Most Probable Paths Using Onsager-Machlup's Method}

As mentioned above, maximal likely trajectories can help us to understand the dynamics of stochastic ice sheet system starting at different initial states. However this doesn't provide any help if we aim for a particular final state. We would like to investigate how we would go from one initial state to a final state of our choice. This is a matter of finding the most probable path for two fixed points. The two fixed points of course have a particular interest as we choose them to be the metastable states. In this subsection, we use Onsager-Machlup's method \cite{OM} to study the most probable transition path connecting initial and final states. From now on we will concentrate exclusively on the transition connecting two metastable states.

Note that system (\ref{Brownian}) is an SDE with multiplicative noise, and the Onsager-Machlup function for SDEs with multiplicative noise can be referred to Bach et al. \cite{Bach1976}. However, it is more convenient to apply Onsager-Machlup's method for an SDE with additive noise numerically \cite{OM}. To this end, we make a transformation: $Z=2\sqrt{X}$ ($X>0$). By It\^o formula, the system (\ref{Brownian}) can be suitably converted to an SDE with additive noise:
\begin{equation}\label{Eqn-Z}
dZ=F(Z)dt+\varepsilon dB_t,
\end{equation}
where $f(Z)=-{\frac{\beta\lambda}{\sqrt{2\sigma}}} \Big(\frac{3}{16}Z^{2}-r \Big)+\frac{1}{6}\beta Z-\frac{\beta^2\varepsilon_0^2}{4\sigma}{Z}^{-1}$ and $\varepsilon= \frac{\beta\varepsilon_0}{\sqrt{2\sigma}}$. Note that, under the variable transformation, $F(Z)=0$ is equivalent to $f(X)-\frac{1}{4}\varepsilon^2=0$. The extra term with respect to $\varepsilon$ appears here. 
It means that equilibrium states for the deterministic counterpart of additive noise system (\ref{Eqn-Z}) are affected by the strength of the noise. As we assume that the noise strength is small, we next look into the transitions between $2\sqrt{X_0}$ and $2\sqrt{X_{+}}$ for system (\ref{Eqn-Z}). By inverse transformation, we then have a better understanding for the transitions between the two metastable states of the original system (\ref{Eqn-L}).\\
\par
For $0\leqslant t_0 \leqslant t_1$, let $\mathcal{T}=\{ z\in C([t_0,t_1];\mathbb{R}):z(t_0)=z_0, z(t_1)=z_1\} $ denote the set of trajectories connecting a point $(t_0,z_0)$ on $2\sqrt{L_0}$ and a point $(t_1,z_1)$ on $2\sqrt{L_{+}}$. Consider an infinitesimally small tubular neighborhood of $z$, 
$$
K(z,\delta)=\{ {z}^{\prime}\in C([t_0,t_1]:|z-z^{\prime}|\leqslant \delta, \text{for }z \in C([t_0,t_1], \delta>0\}.
$$
If $z\in\mathcal{T}$ is differentiable, then the measure of the transition paths lying in the small tubular neighborhood satisfies
\begin{equation}
\mu(K(z,\delta))\propto C(\delta) \int_{\mathcal{T}}\exp\Big(-\frac{1}{2\varepsilon^2}\int_{t_0}^{t_1}OM(z,\dot{z})dt \Big)d\mu_W[z],
\end{equation}
where $C(\delta)$ is a constant with respect to $\delta$, $\mu_W[z]$ is the Wiener measure \cite{OM,Duan}, symbol $\propto$ denotes the proportionality relation, and $OM(z,\dot{z})$ is the so called Onsager-Machlup function which is given by:
\begin{equation}
OM(z,\dot{z})=\big(\dot{z}-F(z)\big)^2+\varepsilon^2F^{\prime}(z).
\end{equation}
We remark that the probability of a transition event lying within a tube along a smooth path $z\in \mathcal{T}$ can be determined by integrating $\mu_\varepsilon[z]$ over this tube. The \emph{most probable path} or \emph{optimal path} connecting points $(t_0,z_0)$ and $(t_1,z_1)$ is defined as minimizers of the OM functional $I_\varepsilon:\mathcal{A}\to \mathbb{R}$ given by
\begin{align}\label{OM functional}
I_\varepsilon[z]:=\int_{t_0}^{t_1}OM(z,\dot{z})dt,
\end{align}
where $\mathcal{A}=\{ z\in H^1([t_0,t_1];\mathbb{R}):z(t_0)=z_0, z(t_1)=z_1\}$. In analogy to classical mechanics, we also call the OM function the Lagrangian function and the OM functional the action functional \cite{DB}. 
\par
Indeed, if we restrict ourselves to twice differentiable functions $z(t)$, the most probable path can be found by variation of the OM functional $I_\varepsilon$. We thus get the Euler-Lagrange equation
\begin{eqnarray}\label{EL}
\frac{d}{dt}\frac{\partial OM(z,\dot{z})}{\partial\dot{z}}=\frac{\partial OM(z,\dot{z})}{\partial z},
\end{eqnarray}
i.e.,
\begin{equation}\label{evolu}
\ddot{z}=F^{\prime}(z)F(z)+\frac{\varepsilon^2}{2}F^{\prime\prime}(z).
\end{equation}
with boundary conditions $z(t_0)=z_0, z(t_1)=z_1$.
Note that the problem of solving (\ref{evolu}) is referred to as the two-point boundary value problem. By Theorem 10 in \cite[Section 8.2.5]{Hasty2000}, if there exists a twice differentiable solution to the Euler-Lagrange equation, then it is indeed a (local) minimizer. As this boundary value problem may not have a solution, the most probable transition pathway may not exist.
\par
 We rewrite equation (\ref{evolu}) in Hamiltonian form
\begin{equation}
  \label{Hamiltonian}
\begin{split}
\dot{z}=&\Phi+F(z), \\
\dot{\Phi}=&-F(z)\Phi+\frac{\varepsilon^2}{2}F^{\prime\prime}(z),
\end{split}
\end{equation}
with Hamiltonian function $H(z,p)=\frac{\Phi^2}{2}+F(z)\Phi-\frac{\varepsilon^2}{2}F^{\prime\prime}(z)$. Here the momentum variable $\Phi=\dot{z}-F(z)$ measures the deviation from the deterministic flow. The corresponding Freidlin-Wentzell functional \cite{Chen2018} can be expressed in terms of $\Phi$ as $I_{F}=\int_{t_0}^{t_1}\Phi^2(t)dt$. Therefore, as $\varepsilon\to 0$, most probable paths are well approximated by heteroclinic orbits of (\ref{Hamiltonian}) connecting deterministic solutions and minimizers of $I_\varepsilon$ converge uniformly to minimizers of $I_{F}$ \cite{Dykman2001,Chen2018}. \\
\par
\begin{figure}
 \subfigure[]{  
  \includegraphics[width=2.15in]{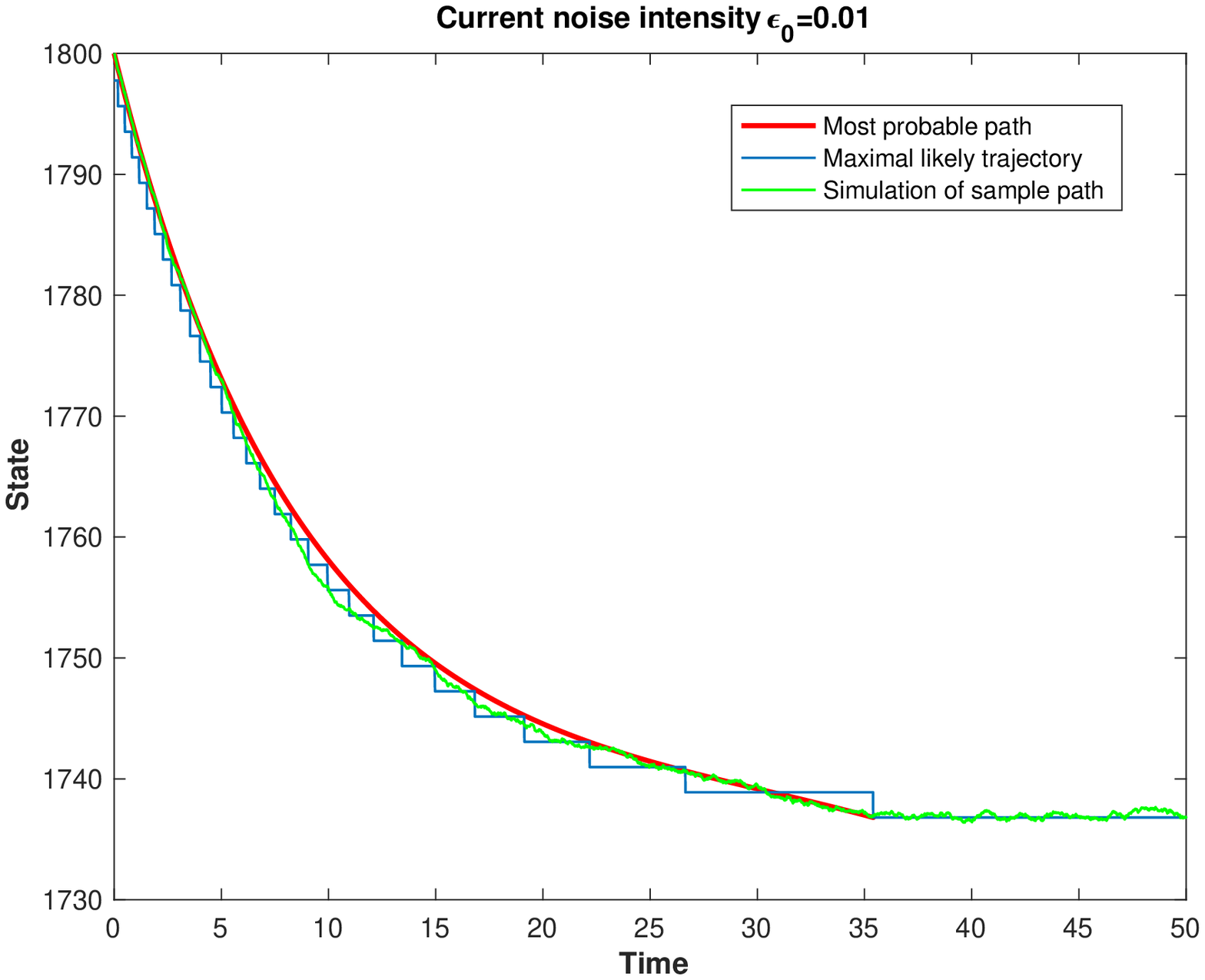}  
  }
   \subfigure[]{  
  \includegraphics[width=2.28in]{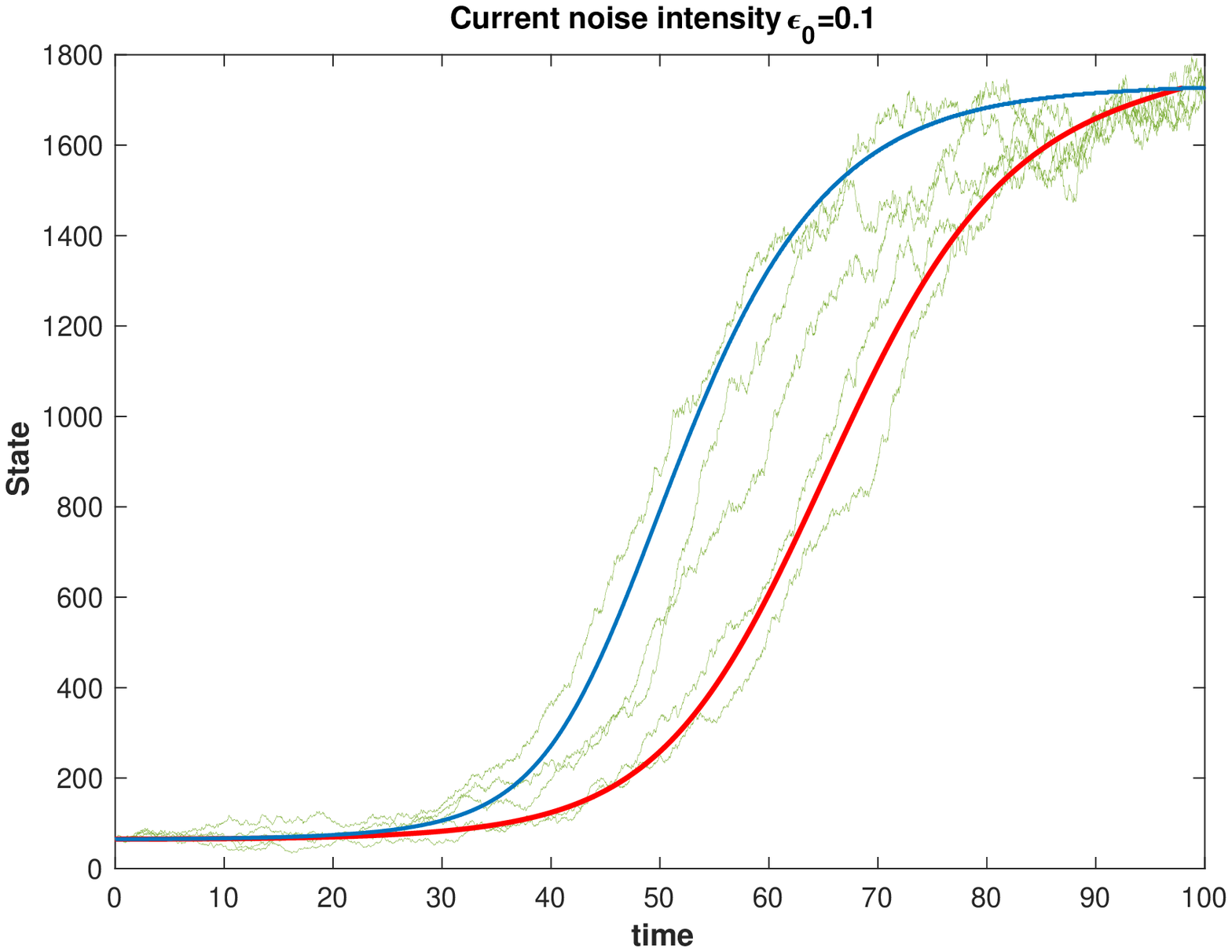}  
  }
	\caption{(Online color) The most probable path (red line) versus the maximal likely trajectory (blue line) for (\ref{Eqn-L}): (a) $\varepsilon_0=0.01$, initial point $X(0)=1800$, final point $X(35.42)=1736.8$; (b) $\varepsilon_0=0.1$, initial point $X(0)=65$, final point $X(97.94)=1734.7$.}  %
	\label{Versus}   %
\end{figure}

In order to illustrate the effectiveness of this method, for a fixed initial point, we choose the final point as that gained by maximal likely trajectory, and choose the transition time as minimal time arriving at the final point. Note that, usually, these two points are not  metastable states, we choose them only for the consideration of methodology. Then, we compare the most probable path (red line) for the stochastic system (\ref{Brownian}) connecting these two fixed points (but not necessarily that they are metastable states here) with the corresponding maximal likely trajectory (blue line) in Fig. \ref{Versus}(a) and (b). Some simulations of actual sample paths (green lines) nearby these two trajectories are also presented in the figures. In Fig. \ref{Versus}(a),  maximal likely trajectory and most probable path coincide with each other. But, in Fig. \ref{Versus}(b), we see that they have obvious difference.
 
\section{Result}
In the following section, we compute the maximal likely evolution trajectories and most probable transition paths, in order to analyze how the ice sheet is formed and how it melts away. The maximal likely trajectories starting from different initial points and most probable transition paths starting from $X_0=0$ and ending at $X_+=1738.6$, are deterministic estimators as time goes on. We will examine maximal likely trajectories and most probable transition paths when system parameters change followed by a brief discussion on general behavior of the system.

\subsection{Maximal Likely Trajectory}
\begin{figure}  
  \subfigure[]{  
  \includegraphics[width=2.3in]{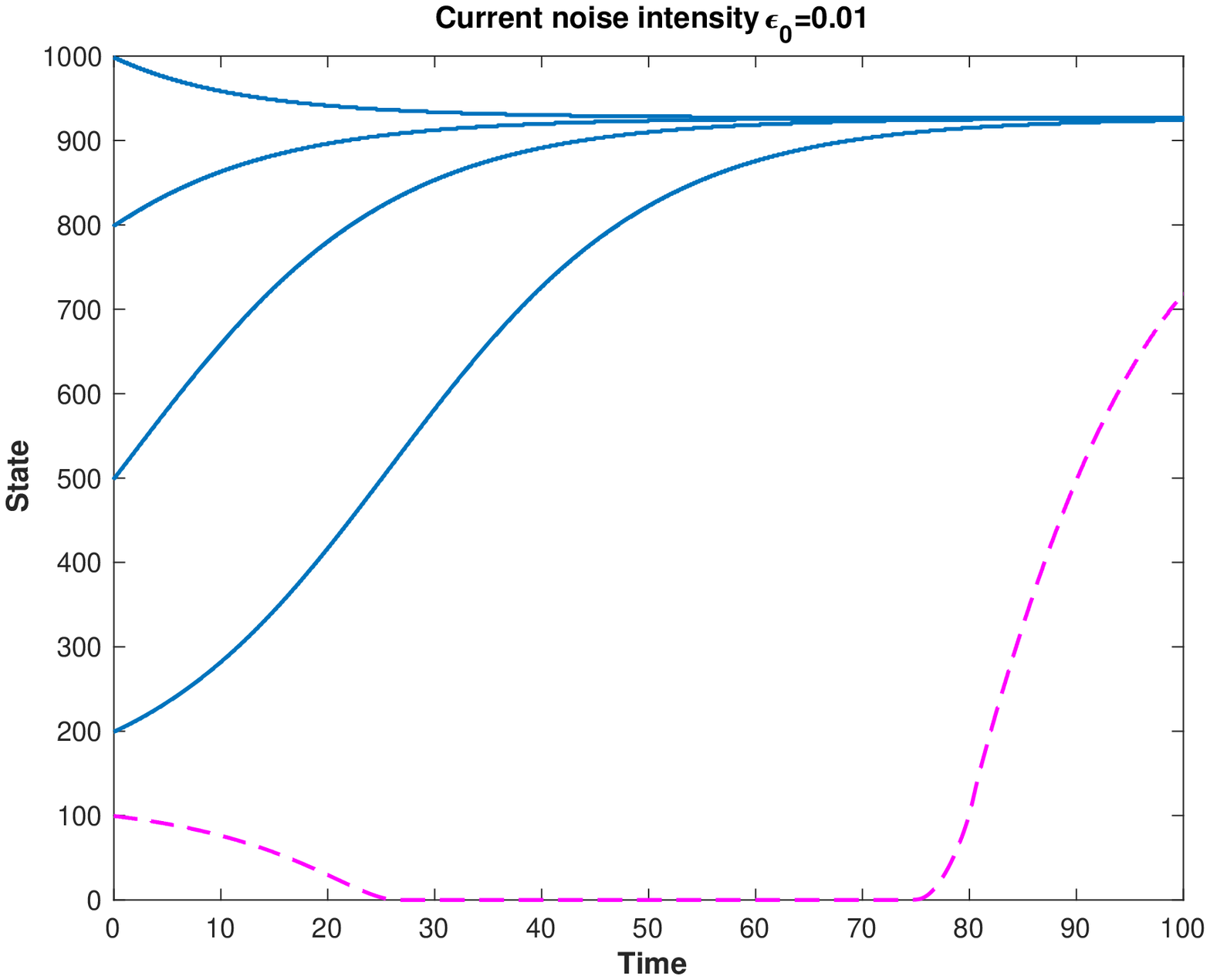}  
  }
  \subfigure[]{  
  \includegraphics[width=2.3in]{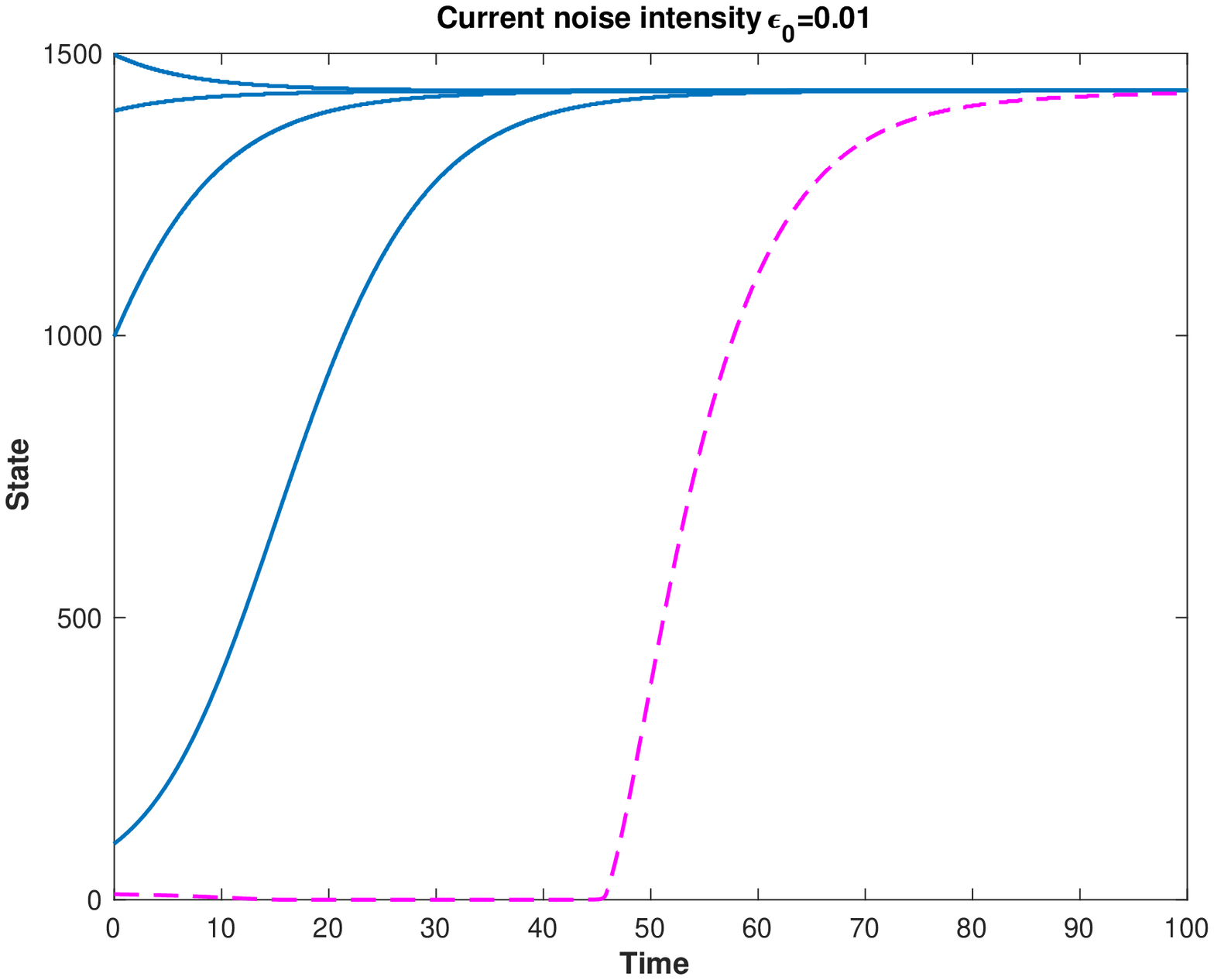}  
  }
  \caption{(Color online) Let $\varepsilon_0=0.01$. Maximal likely evolution trajectories of stochastic system (\ref{Brownian}) with different model parameters: (a) $r=-250(km)$ and $\lambda=0.0012$; (b) $r=-100(km)$ and $\lambda=0.001$.
  }   
  \label{fig:Rlambda} 
\end{figure}

The deterministic dynamical system is bistable in some range of  ice sheet parameter space that contains $\lambda=0.001$, $r=-250$.
For simplicity, we firstly consider five maximum likely trajectories with initial points $X_0=1800$, $1600$, $1000$, $100$ and $50$ to both sides of the equilibria in the deterministic dynamical system. Due to the lengthy computation process time is capped to $T=100$. 
  \begin{figure} 
  \subfigure[]{  
  \includegraphics[width=2.3in]{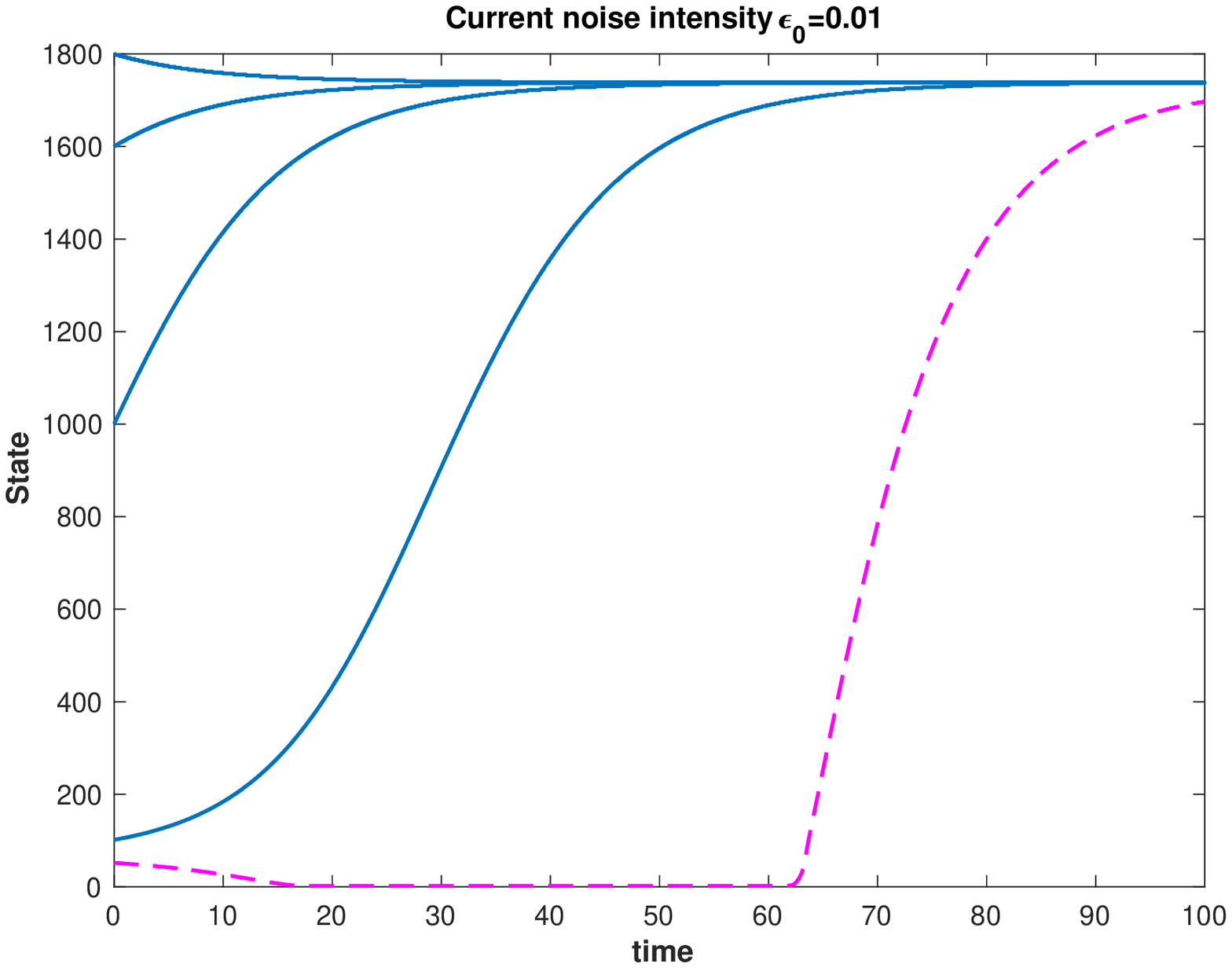}  
  }
  \subfigure[]{  
  \includegraphics[width=2.23in]{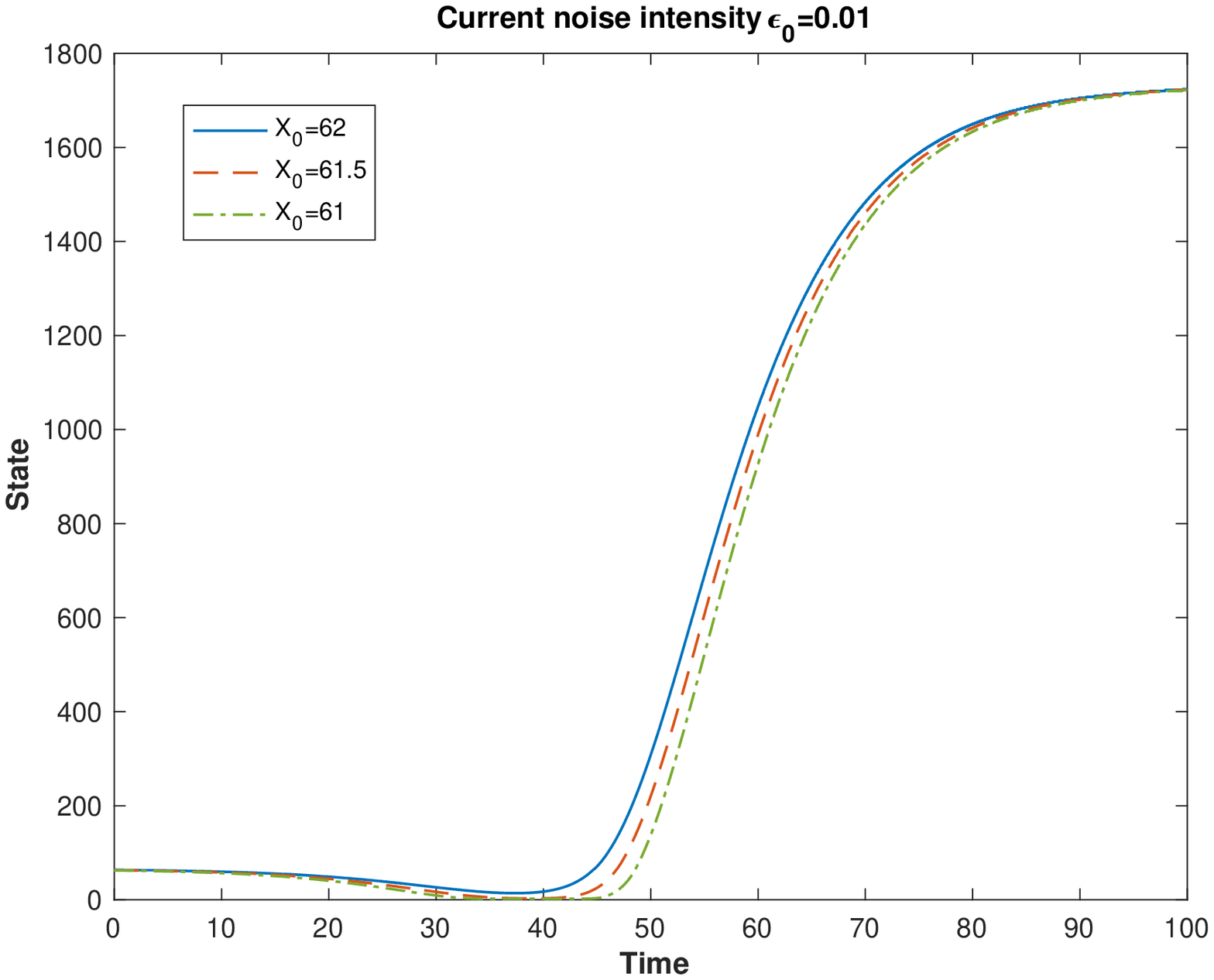}  
  }
  \caption{(Color online) Let $\varepsilon_0=0.01$, $r=-250(km)$ and $\lambda=0.001$.  Maximal likely evolution trajectories of stochastic system (\ref{Brownian}) starting at various initial concentration $X_0$. 
  }   
  \label{fig:MLT} 
\end{figure} 
Fig. \ref{fig:MLT}(a) shows only one maximal likely equilibrium state for stochastic ice sheet system (\ref{Brownian}) with $\varepsilon_0=0.01$, and the value of this maximal likely stable equilibrium state is $1736.8$ which differs slightly from the deterministic stable state $X_+=1738.6$ due to the effect of noise. It is also noticeable that the other deterministic stable state $X_0=0$ is not a maximal likely equilibrium state. We observe that the maximal likely evolution trajectories starting close to the ice-free state (e.g. the pink line in Fig. \ref{fig:MLT}(a)) will go to or go close to zero firstly, but they will reach the ice-covered state ultimately. This could be understood as a reason why the ice sheet concerned here has been formed in its local environment initially. 
As shown in Fig. \ref{fig:MLT}(b), the maximal likely evolution trajectories will reach zero for some time when the values of initial points are smaller than $61.5$.
  \begin{figure} 
  \centering
  \includegraphics[width=2.4in]{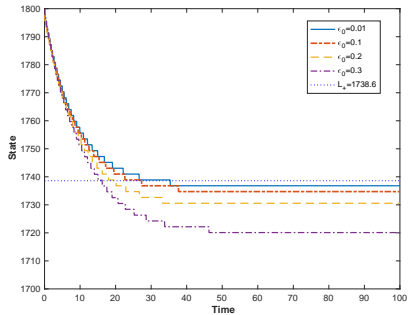}  
  \caption{(Color online) For $\varepsilon_0=0.01$, $r=-250$ and $\lambda=0.001$, the change of maximal likely evolution trajectories of stochastic system (\ref{Brownian}) starting at $X_0=1800$ with respect to different noise intensities.  }   
  \label{fig:diff noise} 
\end{figure} 
Furthermore we investigate the impacts of different levels noise intensity on the maximal likely trajectory for system (\ref{Brownian}) in figure \ref{fig:diff noise}. For convenience, we consider the initial state $X_0=1800$. The maximal likely stable equilibrium state reduces with the increase of the noise intensity as seen in Fig. \ref{fig:diff noise}. Recall that, in Fig. \ref{fig_potential}, the metastable states depend on model  parameters $\lambda$ and $r$. With different model  parameters, the maximal likely stable equilibrium state will also change; see Fig. \ref{fig:Rlambda}.

We end this section with a short summary of key points. First, while the deterministic model has two metastable states, in the stochastic it's reduced to one. In addition this metastable state is dependent both on the parameters and the noise intensity.

\subsection{Most Probable Path}
\begin{figure}
  \subfigure[]{ 
    \includegraphics[width=2.3in]{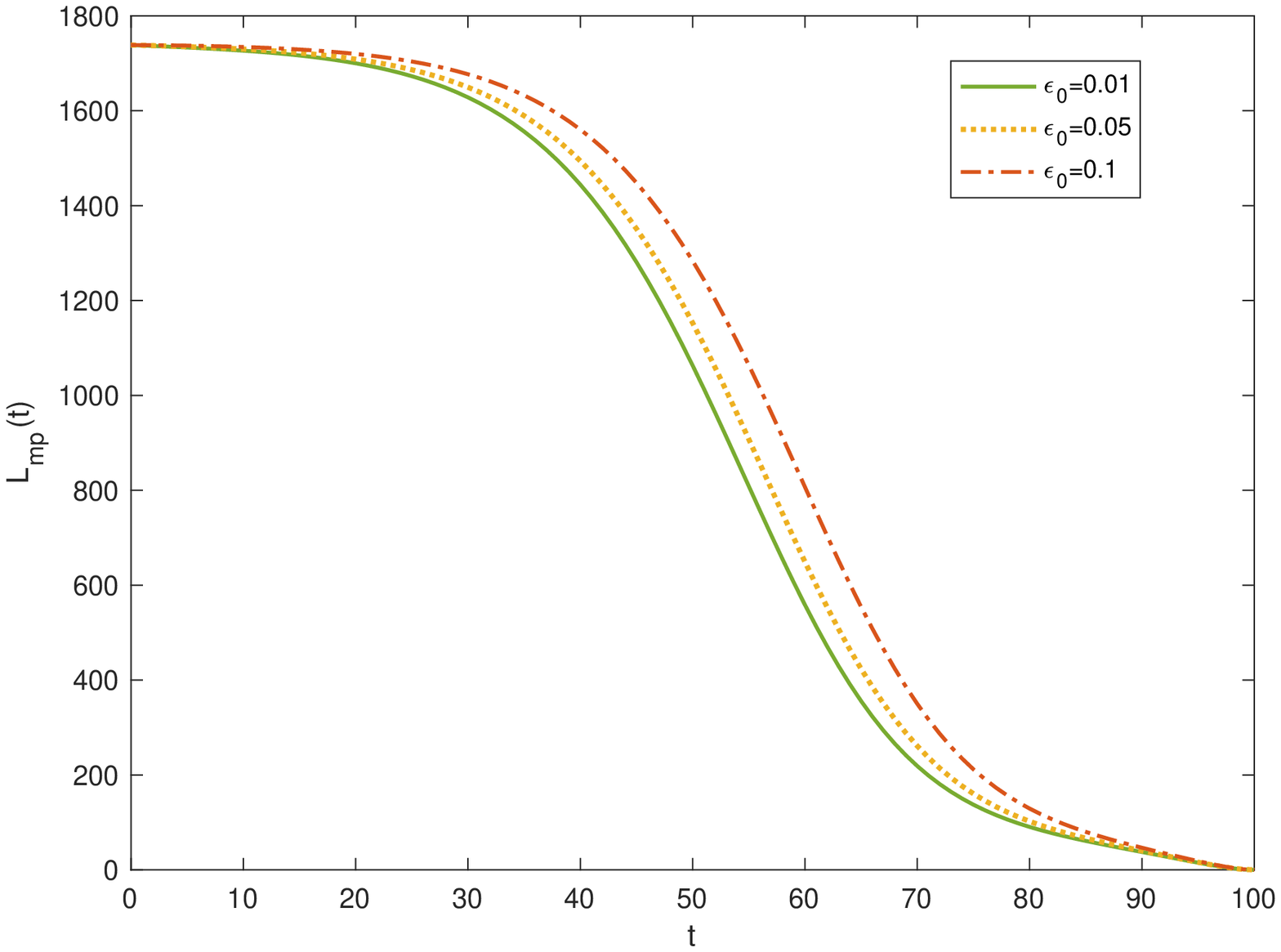}  
  } 
  \subfigure[]{ 
    \includegraphics[width=2.3in]{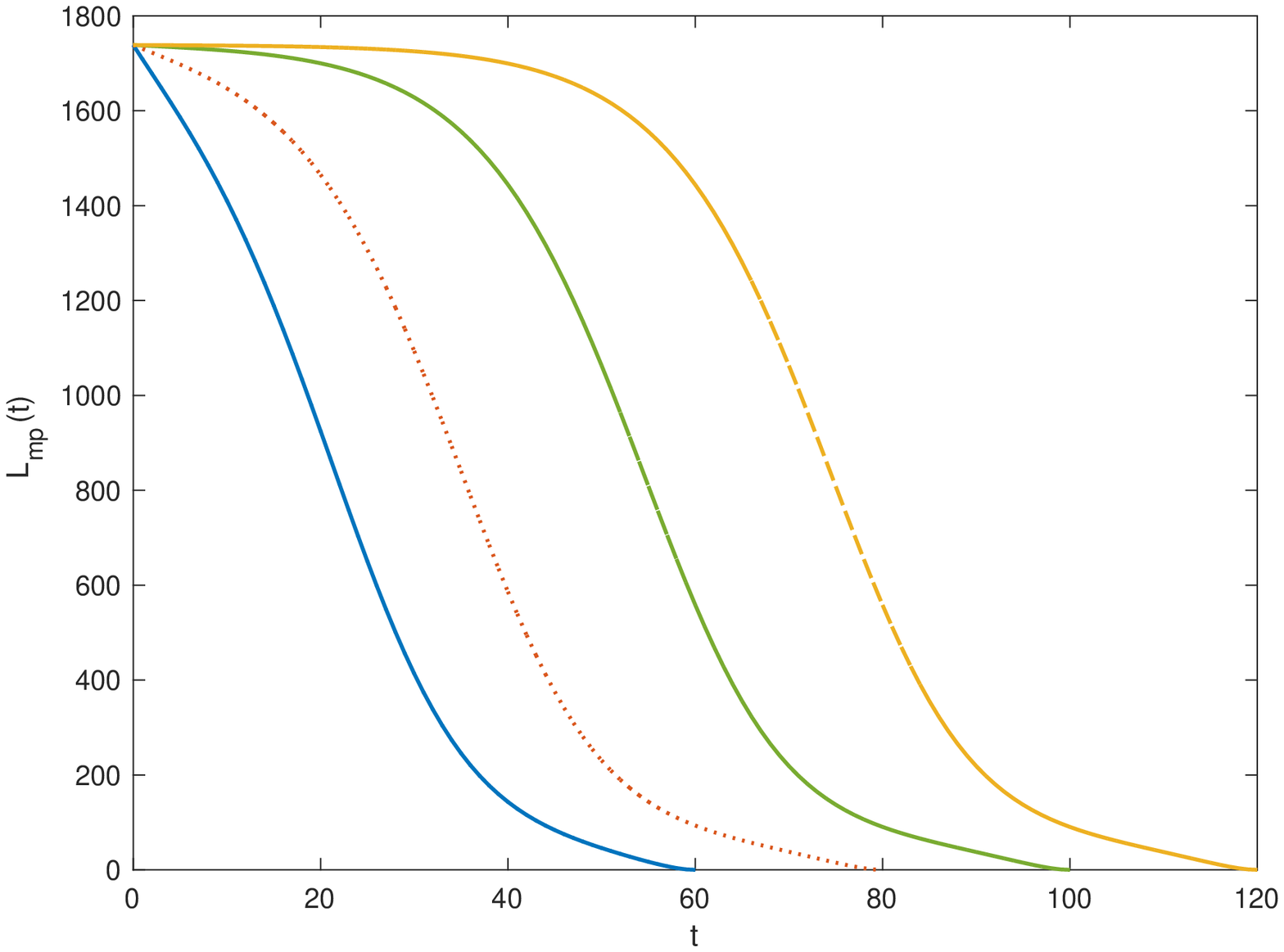} 
  } 
  \caption{(Online color) For system (\ref{Eqn-Z}), set initial value $z(0) =2\sqrt{X_+}$  and final value $z(t_1) =0$. (a) For fixed $t_1=100$, the change of most probable path for original system (\ref{Eqn-L}) with respect to different noise strength $\varepsilon_0$; (b) For fixed $\varepsilon_0=0.01$, the change of most probable path for original system (\ref{Eqn-L}) with respect to different time interval parameters.} 
  \label{MPP} 
\end{figure}

While in the previous section we saw how the trajectory of the system evolves over time, we are also interested to see the trajectory it would follow to reach a specific end.
For stochastic ice sheet system (\ref{Brownian}), we now examine the most probable transition pathway starting at the ice-covered metastable state $X_+=1738.6$ and ending at the ice-free state $0$. As seen in Fig. \ref{MPP}(a), the most probable ice sheet height $L_{mp}$ decreases slowly at first, but after 40-50 kyrs from the start the melting rate accelerates. After that it continues to decrease on a downward trend and passes the `barrier' $X_-=63.9$ after about 80kyrs, and finally, it reaches the ice-free state. The evolution of such most probable transition pathway depends crucially on the noise intensity $\varepsilon_0$ and the system running time $t_1$ as seen in Fig. \ref{MPP}(a) and (b).

Fig. \ref{MPP}(a) shows the most probable path with different noise intensity whilst Fig. \ref{MPP}(b) with different running time. Both affect the most probable path. It should be noted that this shows the most probable path for this transition to happen, not the likeliness for this event. The transition from the ice-covered to the ice-free state would be a rare event. As we have seen in the previous section, Fig. \ref{fig:MLT}, even when starting from very low height the most probable outcome is to reach the metastable state $X_+$.

\section{Conclusion}
 In this work, we have established an ice-sheet model with Brownian noise. It was constructed closely to its deterministic predecessor. We examine the maximal likely trajectories and most probable transition paths, i.e., we visualize the trajectories from different initial state (including ice-free state) to ice-covered state as well as transition pathways from ice-free state to ice covered state. The maximal likely trajectories are calculated via numerically solving the Fokker-Planck equation for the stochastic ice sheet model (\ref{Brownian}). The most probable transition pathways are computed by numerically solving a two-point boundary value problem.
 
 For a stochastic ice-sheet system, we have observed that the maximal likely trajectories starting from near ice-free state would converge to only one maximal likely stable equilibrium state, which could be recognized as a ice-covered state. It shows that there would always be a very thick ice sheet no matter what the initial situation is, and the length of ice sheet would reach definite value (e.g. 1736.8\emph{km} in the case of Fig. \ref{fig:MLT}) after a long time. At the same time, we have also noticed some peculiar or counter-intuitive phenomena. For example, the initial ice cap might melt for some time and then gradually form a big one, if its initial length is small enough (e.g. smaller than 61\emph{km} in the case of Fig. \ref{fig:MLT}). This phenomenon does not occur in the case of deterministic model.
 
 Although it seems that the total melting of the ice-sheet in this model is a very rare event, it is necessary to study the transition from ice-covered state to ice-free state due to this even is more concerned about the people. The method of most probable transition paths by minimizing the Onsager-Machlup action functional could thus be applied. For certain evolution time scale and system parameters, we have indeed observed that the most probable transition pathway exists under Brownian noise. Furthermore, we have characterized the evolution with varying noise parameters. Therefore, we can predict the melting route at a given future time.
 
 The findings in this work may provide helpful insights for further practical research, to verify the change of an ice sheet and so on. And these two methods applied here could be also used in other practical models.
 

\section*{Acknowledgements}
The authors are grateful to Xiaoli Chen, Xiujun Cheng and Yuanfei Huang for helpful discussions and comments. This work was partly supported by the NSF grant 1620449, and NSFC grants 11531006 and  11771449.

\section*{References}

\bibliography{mybibfile}

\end{document}